\documentclass[10pt,a4paper]{article}
\usepackage{CJK}
\usepackage{mathrsfs}
\usepackage{amsfonts}
\begin{document}
\begin{CJK*}{GBK}{song}
\begin{center}
 {\huge\textbf{Some Rapidly Converging Series for $\zeta(2n+1)$ from Abstract Operators}}
 \end{center}
\vspace{0.3cm}

\begin{center}
\rm Guang-Qing Bi
\end{center}

\begin{abstract}
The author derives new family of series representations for the values of the Riemann Zeta function $\zeta(s)$ at positive odd integers. For $n\in\mathbb{N}$, each of these series representing $\zeta(2n+1)$ converges remarkably rapidly with its general term having the order estimate:
$$O(m^{-2k}\cdot k^{-2n+1})\qquad(k\rightarrow\infty;\quad m=3,4,6).$$
The method is based on the mapping relationships between analytic functions and periodic functions using the abstract operators $\cos(h\partial_x)$ and $\sin(h\partial_x)$, including the mapping relationships between power series and trigonometric series, if each coefficient of a power series is respectively equal to that of a trigonometric series. Thus we obtain a general method to find the sum of the Dirichlet series of integer variables. By defining the Zeta function in an abstract operators form, we have further generalized these results on the whole complex plane.
\end{abstract}

\footnotetext{\hspace*{-.45cm} \noindent{\bf 2010 \emph{Mathematics Subject Classification}.}\ \ Primary 11M06, 35S05; Secondary 42A24, 33E20, 40A30.}

\section{Introduction}
\noindent
\renewcommand\theequation{1.\arabic{equation}}

The trigonometric series, especially the Fourier series, is of great importance in both mathematics and physics. Expanding periodic functions into Fourier series has become a very mature theory. During early years of the last century, researchers had recognized the significance of the inverse problem, which is how to find the sum function of a certain Fourier series.

A class of pseudo-differential operators extended from $\mathscr{S}(\mathbb{R}^n)$ to $C^\infty(\Omega)$ is called the abstract operators such as $\cos(h\partial_x)$ and $\sin(h\partial_x)$. The concept of abstract operators is based on the analytic continuity fundamental theorem contained in a 1997 paper entitled \emph{Applications of abstract operators to partial differential equations} (See \cite[p. 7]{bi97}). In other words, if two abstract operators $A$ and $B$ acts on $e^{\xi{x}},\,x\in\mathbb{R}^n,\,\xi\in\mathbb{R}^n$ such that $Ae^{\xi{x}}=Be^{\xi{x}}$, then $Ag(x)=Bg(x),\,\forall{g(x)}\in{C^\infty(\Omega)},\,\Omega\subset\mathbb{R}^n$. Therefore, the abstract operators $f(\partial_x),\,\partial_x:=(\partial_{x_1},\partial_{x_2},\cdots,\partial_{x_n})$ is defined by (See \cite[p. 8]{bi97})
\begin{equation}\label{abs}
 f(\partial_x)e^{\xi{x}}:=f(\xi)e^{\xi{x}},\quad\forall f(\xi)\in C^\infty(\mathbb{R}^n),
\end{equation}
where $f(\xi),\xi\in\mathbb{R}^n$ is called the symbols of abstract operators $f(\partial_x)$.  From this, we can derive the operation rules of abstract operators, such as
\begin{equation}\label{e0}
f(\partial_x)(e^{\xi_0x}g(x))=g(\partial_\xi)(e^{\xi{x}}f(\xi))|_{\xi=\xi_0},
\end{equation}
\begin{equation}\label{e1}
f(\partial_x)(e^{\xi{x}}g(x))=e^{\xi{x}}f(\xi+\partial_x)g(x),
\end{equation}
\begin{equation}\label{e1'}
 e^{h\partial_x}f(x)=f(x+h),\quad h\partial_x:=h_1\partial_{x_1}+h_2\partial_{x_2}+\cdots+h_n\partial_{x_n}
\end{equation}
and so on. The concept of abstract operators and the theory of partial differential equations expressed in terms of abstract operators have been established in \cite{bi18}-\cite{bi11}. We believe that our findings here can be useful for investing the number theoretical properties for $\zeta(2n+1)$. In this paper we describe the operation rules of the abstract operators $\cos(h\partial_x)$ and $\sin(h\partial_x)$, and derive a series of operation formulas for these two operators by their operation rules. In particular, the abstract operators $\cos(h\partial_x)$ and $\sin(h\partial_x)$ can express the mapping relations between periodic functions and analytic functions more explicitly. According to these mapping relations, if we know the sum function of the power series, we can obtain that of the trigonometric series, and the non-analytical points of which are also determined at the same time. Thus we can use the sum of power series to obtain that of corresponding trigonometric series. A summation method of trigonometric series expressed by abstract operators is proposed.

In fact, the trigonometric series has been applied to the Riemann Zeta function by many mathematicians such as Choi \cite{Choi}, Katsurada \cite{Kats} and Tsumura \cite{Tsu} etc. By using the abstract operators $\cos(h\partial_x)$ and $\sin(h\partial_x)$,  we can see that our summation method of trigonometric series is quite suitable to find the sum of the Dirichlet series of integer variables including the Zeta function. The Riemann Zeta function $\zeta(s)$ is defined by
\begin{equation}\label{00}
  \zeta(s):=\sum^\infty_{n=1}\frac{1}{n^s}\qquad(\Re(s)>1),
\end{equation}
which can be continued analytically to the whole $s$-plane except for a simple pole at $s=1$ with its residue $1$.

Recall $\zeta(0)=-1/2$. In 1749, Euler proved
\begin{equation}\label{0}
  \zeta(2n)=(-1)^{n-1}\frac{(2\pi)^{2n}}{2(2n)!}B_{2n}\qquad(n\in\mathbb{N}),
\end{equation}
where $B_m$ are the Bernoulli numbers. Since $B_m$ are rational numbers, the expression (\ref{0}) shows that $\zeta(2n)$ is a rational multiple of $\pi^{2n}$. Yet, there is no analogue for $\zeta(2n+1)$. In 1978, Ap\'{e}ry proved that $\zeta(3)$ is irrational. Since then, it is known that infinitely many $\zeta(2n+1)$ are irrational. Instead, many mathematicians have tried to present rapidly converging series for $\zeta(2n+1)$. In this paper, we aim to derive certain new rapidly converging series for $\zeta(2n+1)$, by using the abstract operators $\cos(h\partial_x)$ and $\sin(h\partial_x)$. We will see that the concept of abstract operators can be applied in defining Riemann Zeta function on the whole complex plane without the experience of analytic continuation process. So we can extend the related results of $\zeta(2n+1)$ to the whole complex plane.

\section{The summation method of Fourier series expressed by abstract operators $\cos(h\partial_x)$ and $\sin(h\partial_x)$}
\noindent\setcounter{equation}{0}
\renewcommand\theequation{2.\arabic{equation}}

Let $h\partial_x=\langle{h,\partial_x}\rangle=h_1\partial_{x_1}+h_2\partial_{x_2}+\cdots+h_n\partial_{x_n}$. Then $\cos(h\partial_x)$ and $\sin(h\partial_x)$ are the abstract operators taking $\cos(hb)$ and $\sin(hb)$ as the symbols respectively, namely (See \cite[p. 8]{bi97})
\begin{equation}\label{1}
    \cos(h\partial_x)e^{bx}:=\cos(bh)e^{bx},\quad\sin(h\partial_x)e^{bx}:=\sin(bh)e^{bx}.
\end{equation}
Here $bx=b_1x_1+b_2x_2+\cdots+b_nx_n,\;bh=b_1h_1+b_2h_2+\cdots+b_nh_n$. Further, their operation rules can be expressed as the following three groups of operator relationships by Guang-Qing Bi \cite[pp. 7-9]{bi97}:

\textbf{Theorem 2.1.} (See \cite[p. 9, Theorem 3, 4, and 6]{bi97}) Let $x\in\mathbb{R}^n,\;h\in\mathbb{R}^n$, $h\partial_x=\langle{h,\partial_x}\rangle=h_1\partial_{x_1}+\cdots+h_n\partial_{x_n}$. For the abstract operators $\cos(h\partial_x)\,\mbox{and}\,\sin(h\partial_x)$, we have
\begin{equation}\label{y0}
  \cos(h\partial_x)f(x)=\Re[f(x+ih)],\quad\sin(h\partial_x)f(x)=\Im[f(x+ih)],
\end{equation}
$\forall{f(z)}\in{C}^\infty(\Omega),\,z=x+iy\in\Omega\subset\mathbb{C}^n$;

\parbox{11cm}{\begin{eqnarray*}\label{y1}
                \sin(h\partial_x)(uv) &=& \cos(h\partial_x)v\cdot\sin(h\partial_x)u+\sin(h\partial_x)v\cdot\cos(h\partial_x)u,\\
                \cos(h\partial_x)(uv) &=& \cos(h\partial_x)v\cdot\cos(h\partial_x)u-\sin(h\partial_x)v\cdot\sin(h\partial_x)u;
              \end{eqnarray*}}\hfill\parbox{1cm}{\begin{eqnarray}\end{eqnarray}}

\parbox{11cm}{\begin{eqnarray*}\label{y2}
                \sin(h\partial_x)\frac{u}{v} &=& \frac{\cos(h\partial_x)v\cdot\sin(h\partial_x)u-\sin(h\partial_x)v\cdot\cos(h\partial_x)u}
                {(\cos(h\partial_x)v)^2+(\sin(h\partial_x)v)^2}, \\
                \cos(h\partial_x)\frac{u}{v} &=& \frac{\cos(h\partial_x)v\cdot\cos(h\partial_x)u+\sin(h\partial_x)v\cdot\sin(h\partial_x)u}
                {(\cos(h\partial_x)v)^2+(\sin(h\partial_x)v)^2}.
              \end{eqnarray*}}\hfill\parbox{1cm}{\begin{eqnarray}\end{eqnarray}}

\textbf{Theorem 2.2.} (See \cite[p. 9, Theorem 5]{bi97}) Let $h_0\in\mathbb{R},\;x(t)\in\mathbb{R}^n,\;t\in\mathbb{R}^1,\;X\in\mathbb{R}^n,\,Y\in\mathbb{R}^n$,
$Y\partial_X=\langle{Y,\partial_X}\rangle=Y_1\partial_{X_1}+\cdots+Y_n\partial_{X_n}$. Then we have

\parbox{11cm}{\begin{eqnarray*}\label{y3}
 \sin(h_0\partial_t)f(x(t)) &=& \sin(Y\partial_X)f(X),\\
 \cos(h_0\partial_t)f(x(t)) &=& \cos(Y\partial_X)f(X),
\end{eqnarray*}}\hfill\parbox{1cm}{\begin{eqnarray}\end{eqnarray}}
where $X_j=\cos(h_0\partial_t)x_j(t),\;Y_j=\sin(h_0\partial_t)x_j(t),\;j=1,\cdots,n$.

In the special case when $n=1$, (\ref{y3}) can easily be restated as

\parbox{11cm}{\begin{eqnarray*}\label{yb1}
\sin\left(h_0\frac{d}{dt}\right)f(x(t)) &=& \sin\left(Y\frac{\partial}{\partial{X}}\right)f(X),\\
\cos\left(h_0\frac{d}{dt}\right)f(x(t)) &=& \cos\left(Y\frac{\partial}{\partial{X}}\right)f(X),
\end{eqnarray*}}\hfill\parbox{1cm}{\begin{eqnarray}\end{eqnarray}}
where
$Y=\sin(h_0\partial_t)x(t),\;X=\cos(h_0\partial_t)x(t),\;t\in\mathbb{R}^1,\;h_0\in\mathbb{R}$.

\textbf{Theorem 2.3.} (See \cite[p. 9, Theorem 7]{bi97}) Let $u=g(y)$ be a monotonic function on its domain. If $y=f(bx)$ is the inverse function of $bx=g(y)$ such that $g(f(bx))=bx$, where $bx=b_1x_1+b_2x_2+\cdots+b_nx_n,\;bh=b_1h_1+b_2h_2+\cdots+b_nh_n$, then $\sin(h\partial_x)f(bx)$ (denoted by $Y$) and
$\cos(h\partial_x)f(bx)$ (denoted by $X$) can be determined by the following set of equations:
\begin{equation}\label{y4}
    \left\{\begin{array}{l@{\qquad}l}\displaystyle
    \cos\left(Y\frac{\partial}{\partial{X}}\right)g(X)=bx,&x\in\mathbb{R}^n,\;b\in\mathbb{R}^n,\\\displaystyle
    \sin\left(Y\frac{\partial}{\partial{X}}\right)g(X)=bh, &h\in\mathbb{R}^n.
    \end{array}\right.
\end{equation}

When acting on elementary functions, the abstract operators $\cos(h\partial_x)$ and $\sin(h\partial_x)$ have complete operation formulas as follows:

By using (\ref{y0})  we have

\parbox{10cm}{\begin{eqnarray*}\label{2a}
                \cos(h\partial_x)\sin{bx} &=& \cosh(bh)\sin{bx}, \\
                \sin(h\partial_x)\sin{bx} &=& \sinh(bh)\cos{bx}.
              \end{eqnarray*}}\hfill\parbox{1cm}{\begin{eqnarray}\end{eqnarray}}

\parbox{10cm}{\begin{eqnarray*}\label{2b}
                \cos(h\partial_x)\cos{bx} &=& \cosh(bh)\cos{bx}, \\
                \sin(h\partial_x)\cos{bx} &=& -\sinh(bh)\sin{bx}.
              \end{eqnarray*}}\hfill\parbox{1cm}{\begin{eqnarray}\end{eqnarray}}

Based on (\ref{2a}) and (\ref{2b}), by using (\ref{y2}) we obtain

\parbox{10cm}{\begin{eqnarray*}\label{2c}
                \cos(h\partial_x)\tan{bx} &=& \frac{\sin(2bx)}{\cosh(2bh)+\cos(2bx)}, \\
                \sin(h\partial_x)\tan{bx} &=& \frac{\sinh(2bh)}{\cosh(2bh)+\cos(2bx)}.
              \end{eqnarray*}}\hfill\parbox{1cm}{\begin{eqnarray}\end{eqnarray}}

\parbox{10cm}{\begin{eqnarray*}\label{4'}
                \cos(h\partial_x)\cot{bx} &=& \frac{\sin(2bx)}{\cosh(2bh)-\cos(2bx)},\\
                \sin(h\partial_x)\cot{bx} &=& \frac{\sinh(2bh)}{\cos(2bx)-\cosh(2bh)}.
              \end{eqnarray*}}\hfill\parbox{1cm}{\begin{eqnarray}\end{eqnarray}}

By using (\ref{1}) and (\ref{y4}), we have
\begin{eqnarray*}
    \left\{\begin{array}{l@{\qquad}l}\displaystyle
    e^X\cos{Y} = bx, & X=\cos(h\partial_x)\ln(bx),\\\displaystyle
    e^X\sin{Y} = bh, & Y=\sin(h\partial_x)\ln(bx).
    \end{array}\right.
\end{eqnarray*}
By solving this set of equations we obtain

\parbox{10cm}{\begin{eqnarray*}\label{7}
                \cos(h\partial_x)\ln(bx) &=& \frac{1}{2}\ln\left((bx)^2+(bh)^2\right), \\
                \sin(h\partial_x)\ln(bx) &=&\textrm{arccot}\frac{bx}{bh}.
              \end{eqnarray*}}\hfill\parbox{1cm}{\begin{eqnarray}\end{eqnarray}}

By using (\ref{2c}) and (\ref{y4}), we have
\begin{eqnarray*}
    \left\{\begin{array}{l@{\qquad}l}\displaystyle
    \frac{\sin2X}{\cosh2Y+\cos2X}=bx, & X=\cos(h\partial_x)\arctan{bx},\\\displaystyle
    \frac{\sinh2Y}{\cosh2Y+\cos2X}=bh, & Y=\sin(h\partial_x)\arctan{bx}.
    \end{array}\right.
\end{eqnarray*}
By solving this set of equations we have
\[1+(bx)^2+(bh)^2=1+\frac{\sin^22X+\sinh^22Y}{(\cosh2Y+\cos2X)^2}=\frac{2\cosh2Y}{\cosh2Y+\cos2X}=\frac{2bh\cosh2Y}{\sinh2Y}.\]
and
\[1-(bx)^2-(bh)^2=1-\frac{\sin^22X+\sinh^22Y}{(\cosh2Y+\cos2X)^2}=\frac{2\cos2X}{\cosh2Y+\cos2X}=\frac{2bx\cos2X}{\sin2X}.\]
So we have
\[\tanh2Y=\frac{2bh}{1+(bx)^2+(bh)^2}\quad\mbox{and}\quad\tan2X=\frac{2bx}{1-(bx)^2-(bh)^2}.\]
Thus we obtain

\parbox{10cm}{\begin{eqnarray*}\label{8}
                \sin(h\partial_x)\arctan{bx} &=& \frac{1}{2}\textrm{tanh}^{-1}\frac{2bh}{1+(bx)^2+(bh)^2}, \\
                \cos(h\partial_x)\arctan{bx} &=& \frac{1}{2}\arctan\frac{2bx}{1-(bx)^2-(bh)^2}.
              \end{eqnarray*}}\hfill\parbox{1cm}{\begin{eqnarray}\end{eqnarray}}

Similarly,

\parbox{10cm}{\begin{eqnarray*}\label{9}
                \sin(h\partial_x)\textrm{arccot}\,bx &=& -\frac{1}{2}\textrm{coth}^{-1}\frac{1+(bx)^2+(bh)^2}{2bh}, \\
                \cos(h\partial_x)\textrm{arccot}\,bx &=& \frac{1}{2}\textrm{arccot}\frac{(bx)^2+(bh)^2-1}{2bx}.
              \end{eqnarray*}}\hfill\parbox{1cm}{\begin{eqnarray}\end{eqnarray}}

\parbox{10cm}{\begin{eqnarray*}\label{11}
                \sin(h\partial_x)\sqrt{bx} &=& \sqrt{\frac{\sqrt{(bx)^2+(bh)^2}-bx}{2}}, \\
                \cos(h\partial_x)\sqrt{bx} &=&\sqrt{\frac{\sqrt{(bx)^2+(bh)^2}+bx}{2}}.
              \end{eqnarray*}}\hfill\parbox{1cm}{\begin{eqnarray}\end{eqnarray}}

For hyperbolic and inverse hyperbolic functions, we can also derive the corresponding basic formulas. For instance, correspondingly to
(\ref{2a}) and (\ref{2b}), we have

\parbox{10cm}{\begin{eqnarray*}\label{10a}
                \cos(h\partial_x)\sinh{bx} &=& \cos(bh)\sinh{bx}, \\
               \sin(h\partial_x)\sinh{bx} &=& \sin(bh)\cosh{bx}.
              \end{eqnarray*}}\hfill\parbox{1cm}{\begin{eqnarray}\end{eqnarray}}

\parbox{10cm}{\begin{eqnarray*}\label{10b}
                \cos(h\partial_x)\cosh{bx} &=& \cos(bh)\cosh{bx}, \\
                \sin(h\partial_x)\cosh{bx} &=& \sin(bh)\sinh{bx}.
              \end{eqnarray*}}\hfill\parbox{1cm}{\begin{eqnarray}\end{eqnarray}}

\textbf{Theorem 2.4.}  Let $f(x)\in{L^2}[-c,c]$ be the sum function of the Fourier cosine series, and
$g(x)\in{L^2}[-c,c]$ be that of the corresponding Fourier  sine series, namely
\[ f(x)=\sum^\infty_{n=0}a_n\cos\frac{n\pi{x}}{c}\quad\mbox{and}\quad g(x)=\sum^\infty_{n=0}a_n\sin\frac{n\pi{x}}{c},\]
where $c>0,\;x\in\Omega\subset\mathbb{R}^1$. If $S(t)$ is the sum function of the corresponding power series $\sum^\infty_{n=0}a_nt^n$, namely
\[S(t)=\sum^\infty_{n=0}a_nt^n,\quad(t\in\mathbb{R}^1,\;|t|<r,\;0<r<+\infty),\]
then for $a<x<b$, we have the following mapping relationships:

\parbox{11.5cm}{\begin{eqnarray*}\label{12}
                f(x) &=& \left.\cos\left(\frac{\pi{x}}{c}\frac{\partial}{\partial{z}}\right)S(e^z)\right|_{z=0}, \\
                g(x) &=& \left.\sin\left(\frac{\pi{x}}{c}\frac{\partial}{\partial{z}}\right)S(e^z)\right|_{z=0}.
              \end{eqnarray*}}\hfill\parbox{1cm}{\begin{eqnarray}\end{eqnarray}}
Here interval $(a,b)\subset\Omega$. The endpoints $a$ and $b$ of the interval $a<x<b$ are non-analytical points (singularities) of Fourier series, which can
be uniquely determined by the detailed computation of the right-hand side of (\ref{12}).

\textbf{Proof.} By substituting $S(e^z)=\sum^\infty_{n=0}a_ne^{nz}$ into (\ref{12}), we can prove Theorem 2.4 easily.

\textbf{Note 2.1.} It is easily seen from (\ref{12}) that
\begin{equation}\label{12ct1}
 \left.\cos\left(\frac{\pi{x}}{c}\frac{\partial}{\partial{z}}\right)S(e^z)\right|_{z=0}=\sum^\infty_{n=0}a_n\cos\frac{n\pi{x}}{c};
\end{equation}
\begin{equation}\label{12ct2}
 \left.\sin\left(\frac{\pi{x}}{c}\frac{\partial}{\partial{z}}\right)S(e^z)\right|_{z=0}=\sum^\infty_{n=0}a_n\sin\frac{n\pi{x}}{c}.
\end{equation}
Here $\Omega\in\mathbb{R}^1$. The $\Omega$ can be uniquely determined by the detailed computation of the left-hand side of (\ref{12ct1}) and (\ref{12ct2}) respectively.
For example, according to the proof of Lemma 3.2 in this paper, if $S(e^z)=-\ln(1-e^z)$, then we have $\Omega:=\{x\in\mathbb{R}^1|\sin(\pi{x}/(2c))>0\}$ for (\ref{12ct1}) and $\Omega:=\{x\in\mathbb{R}^1|\sin(\pi{x}/(2c))\neq0\}$ for (\ref{12ct2}) respectively.

In the application of Theorem 2.4, the following theorem can be particularly useful:

\textbf{Theorem 2.5.}  Let $S(t)$ be an arbitrary analytic function integrable in the interval $[0,1]$. Then we have
\begin{equation}\label{13}
  \left.\cos\left(\frac{\pi{x}}{c}\frac{\partial}{\partial{z}}\right)\int^{e^z}_0\!\!\!S(e^z)\,de^z\right|_{z=0}
 = \int^1_0\!\!S(\xi)\,d\xi-\frac{\pi}{c}\int^x_0\!\!\left.\sin\left(\frac{\pi{x}}{c}\frac{\partial}{\partial{z}}\right)
 [S(e^z)\,e^z]\right|_{z=0}dx;
\end{equation}
\begin{equation}\label{14}
 \left.\sin\left(\frac{\pi{x}}{c}\frac{\partial}{\partial{z}}\right)\int^{e^z}_0\!\!S(e^z)\,de^z\right|_{z=0}=
\frac{\pi}{c}\int^x_0\!\!\left.\cos\left(\frac{\pi{x}}{c}\frac{\partial}{\partial{z}}\right)[S(e^z)\,e^z]\right|_{z=0}dx.
\end{equation}

\textbf{Proof.} According to the analytic continuous fundamental theorem, we only need to prove this set of formulas when $S(t)=t^n$, $n\in\mathbb{N}_0:=\{0,1,2,\cdots\}$, this is obvious.

\textbf{Note 2.2.} The analytic continuous fundamental theorem is contained in a 1997 paper entitled \emph{Applications of abstract operators to partial differential equations} by Guang-Qing Bi \cite[p. 7]{bi97}.

\textbf{Theorem 2.6.}   Let $f(x)\in{L^2}[-c,c]$ be the sum function of the Fourier cosine series, and $g(x)\in{L^2}[-c,c]$ be that of the corresponding Fourier sine series, namely
\[f(x)=\sum^\infty_{n=1}a_n\cos\frac{n\pi{x}}{c}\quad\mbox{and}\quad g(x)=\sum^\infty_{n=1}a_n\sin\frac{n\pi{x}}{c},\]
where $c>0,\;x\in\Omega\subset\mathbb{R}^1$. If $S(t)$ is an analytic function integrable in the interval $[0,1]$, and
$\int^t_0S(t)\,dt$ is the sum function of the corresponding power series $\sum^\infty_{n=1}a_nt^n$, namely
\[\int^t_0\!\!S(t)\,dt=\sum^\infty_{n=1}a_nt^n\qquad(t\in\mathbb{R}^1,\;|t|<r,\;0<r<+\infty),\]
then for $a<x<b$ we have the following mapping relationships:

\parbox{11.5cm}{\begin{eqnarray*}\label{15}
f(x) &=&
\int^1_0\!\!S(\xi)\,d\xi-\frac{\pi}{c}\int^x_0\!\!\left.\sin\left(\frac{\pi{x}}{c}\frac{\partial}{\partial{z}}\right)
[S(e^z)\,e^z]\right|_{z=0}dx, \\
g(x)&=&\frac{\pi}{c}\int^x_0\!\!\left.\cos\left(\frac{\pi{x}}{c}\frac{\partial}{\partial{z}}\right)[S(e^z)\,e^z]\right|_{z=0}dx.
\end{eqnarray*}}\hfill\parbox{1cm}{\begin{eqnarray}\end{eqnarray}}
Here interval $(a,b)\subset\Omega$. The endpoints $a$ and $b$ of the interval $a<x<b$ are non-analytical points (singularities) of Fourier series, which can
be uniquely determined by the detailed computation of the right-hand side of (\ref{15}).

\textbf{Proof.} Combining Theorem 2.4 with Theorem 2.5 will lead us to the proof.

\textbf{Example.}   Let $\omega$ be a positive real number, $f(\omega t)\in{L^2}[-\pi,\pi]$ defined by
$$f(\omega t):=\sum^\infty_{n=1}\frac{(-1)^{n-1}}{(3n-1)(3n+1)}\cos(3n\omega{t})\quad\omega{t}\in\Omega\subset\mathbb{R}^1.$$
Then $f(\omega t)=(\pi/(3\sqrt{3}))\cos\omega{t}-1/2$ for $\omega{t}\in(-\pi/3,\pi/3)\subset\Omega$, namely
\begin{equation}\label{22}
    \cos\omega{t}=\frac{3\sqrt{3}}{\pi}\left[\frac{1}{2}+\sum^\infty_{n=1}\frac{(-1)^{n-1}}{(3n-1)(3n+1)}\cos(3n\omega{t})\right]\qquad
    -\frac{\pi}{3}<\omega{t}<\frac{\pi}{3}.
\end{equation}
Here $\Omega:=\{\omega{t}\in\mathbb{R}^1|2\cos\omega{t}-1\neq0\;\mbox{and}\;\cos(\omega{t}/2)\neq0\}$.

\textbf{Proof.}  By using Theorem 2.4, let
\begin{eqnarray*}
  S(x) &=& \sum^\infty_{n=1}(-1)^{n-1}\frac{x^{3n}}{(3n-1)(3n+1)}\qquad{x}\in\mathbb{R}^1,\;|x|\leq1,\\
  f(\omega t) &=& \sum^\infty_{n=1}\frac{(-1)^{n-1}}{(3n-1)(3n+1)}\cos(3n\omega{t})\qquad\omega{t}\in\Omega\subset\mathbb{R}^1.
\end{eqnarray*}
As it is difficult to obtain the sum function $S(x)$ directly, we can use the operators $\frac{d}{dx}(x\cdot)$ and
$\frac{d}{dx}(\frac{1}{x}\cdot)$ to transform the power series into the series familiar to us, namely
\[\frac{d}{dx}\left(\frac{1}{x}\frac{d}{dx}(xS(x))\right)=\sum^\infty_{n=1}(-1)^{n-1}x^{3n-2}=\frac{x}{1+x^3}.\]
Then we can obtain $S(x)$, and then $f(\omega t)$. To describe concisely, the following results will be given directly:
\begin{eqnarray*}
  S(x) &=& \frac{1}{x}\int^x_0\!\!x\,dx\int^x_0\!\!\frac{x}{1+x^3}\,dx=\left(\frac{x}{12}-\frac{1}{12x}\right)\ln\left(x^2-x+1\right)\\
       & & -\left(\frac{x}{6}-\frac{1}{6x}\right)\ln(1+x)+
       \left(\frac{x}{4}+\frac{1}{4x}\right)\frac{2}{\sqrt{3}}\left(\arctan\frac{2x-1}{\sqrt{3}}+\frac{\pi}{6}\right)-\frac{1}{2}.
\end{eqnarray*}

Applying (\ref{1}), (\ref{y1}), (\ref{yb1}), (\ref{7}), (\ref{8}), (\ref{10a}) and (\ref{10b}), we have
\begin{eqnarray*}
f(\omega t) &=& \left.\cos\left(\omega{t}\frac{\partial}{\partial{z}}\right)S(e^z)\right|_{z=0} \\
&=& \frac{1}{6}\left.\cos\left(\omega{t}\frac{\partial}{\partial{z}}\right)\left[\,\sinh{z}\ln\left(e^{2z}-e^z+1\right)\right]\right|_{z=0}\\
& &-\frac{1}{3}\left.\cos\left(\omega{t}\frac{\partial}{\partial{z}}\right)[\,\sinh{z}\ln\left(1+e^z\right)]\right|_{z=0} \\
& & +\frac{1}{\sqrt{3}}\left.\cos\left(\omega{t}\frac{\partial}{\partial{z}}\right)\left[\,\cosh{z}\left(\arctan\frac{2e^z-1}{\sqrt{3}}
  +\frac{\pi}{6}\right)\right]\right|_{z=0}-\frac{1}{2}\\
&=& -\frac{1}{6}\sin\omega{t}\,\textrm{arccot}\frac{(2\cos\omega{t}-1)\cos\omega{t}}{(2\cos\omega{t}-1)\sin\omega{t}}
+\frac{1}{3}\sin\omega{t}\,\textrm{arccot}\frac{\cos^2(\omega{t}/2)}{\sin(\omega{t}/2)\cos(\omega{t}/2)} \\
& & +\frac{1}{2\sqrt{3}}\cos\omega{t}\arctan\frac{\sqrt{3}(2\cos\omega{t}-1)}{2\cos\omega{t}-1}+\frac{\pi}{6\sqrt{3}}\cos\omega{t}-\frac{1}{2}
\qquad\omega{t}\in\Omega.
\end{eqnarray*}

When $2\cos\omega{t}-1\neq0$ and $\cos(\omega{t}/2)\neq0$, we obtain
\begin{eqnarray*}
  f(\omega t) &=& -\frac{1}{6}\sin\omega{t}\,\textrm{arccot}\cot\omega{t}+\frac{1}{3}\sin\omega{t}\,\textrm{arccot}\cot\frac{\omega{t}}{2} \\
       & & +\frac{1}{2\sqrt{3}}\cos\omega{t}\arctan\sqrt{3}+\frac{\pi}{6\sqrt{3}}\cos\omega{t}-\frac{1}{2}
       =\frac{\pi}{3\sqrt{3}}\cos\omega{t}-\frac{1}{2}.
\end{eqnarray*}
Since $2\cos\omega{t}-1\neq0$ and $\cos(\omega{t}/2)\neq0$, there are four non-analytical points: $\omega{t}=-\pi,-\pi/3,\pi/3,\pi$ in the interval $[-\pi,\pi]$. Thus we obtain (\ref{22}) for $-\pi/3<\omega{t}<\pi/3$.

\textbf{Note 2.3.}  For $f(\omega t)\in{L^2}[-\pi,\pi]$, since $2\cos\omega{t}-1\neq0$ and $\cos(\omega{t}/2)\neq0$, there are countable infinite non-analytical points on $-\infty<\omega{t}<\infty$, namely $\Omega:=\{\omega{t}\in\mathbb{R}^1|2\cos\omega{t}-1\neq0\;\mbox{and}\;\cos(\omega{t}/2)\neq0\}$. The open interval $(-\pi/3,\pi/3)\subset\Omega$.

\section{Preliminary applications}
\noindent\setcounter{equation}{0}
\renewcommand\theequation{3.\arabic{equation}}

Let $S_0(t)$ be a function analytic in the neighborhood of $t=0$ and
\[S_0(t)=\sum^\infty_{n=1}a_nt^n\qquad(t\in\mathbb{R}^1,\;|t|<r,\;0<r<+\infty),\]
where $a_n$ are rational numbers. Then the sum function $S_m(t)$ is defined as
\begin{equation}\label{23}
S_m(t):=\underbrace{\int^t_0\frac{dt}{t}\cdots}_m\int^t_0S_0(t)\,\frac{dt}{t}=\sum^\infty_{n=1}a_n\frac{t^n}{n^m}\qquad m\in\mathbb{N}_0:=\{0,1,2,\cdots\}.
\end{equation}
Apparently $S_m(1)$ is the sum function of the Dirichlet series taking $m$ as the variable.

\textbf{Lemma 3.1.} $S_m(t)$ satisfies the following recurrence relation:
\begin{equation}\label{24}
\int^t_0S_{m-1}(t)\,\frac{dt}{t}=S_m(t)\qquad m\in\mathbb{N}:=\{1,2,\cdots\}.
\end{equation}

\textbf{Lemma 3.2.} Let $x\in\mathbb{R}^1$ with $0<x<2c$. Then
\begin{equation}\label{31a}
\left.\sin\left(\frac{\pi{x}}{c}\frac{\partial}{\partial{z}}\right)(-\ln(1-e^z))\right|_{z=0}=
 \frac{\pi}{2}-\frac{\pi{x}}{2c};
\end{equation}
\begin{equation}\label{31b}
 \left.\cos\left(\frac{\pi{x}}{c}\frac{\partial}{\partial{z}}\right)(-\ln(1-e^z))\right|_{z=0}=
 -\ln\left(2\sin\frac{\pi{x}}{2c}\right).
\end{equation}

\textbf{Proof.} By using (\ref{1}), (\ref{yb1}) and (\ref{7}), we have
\begin{eqnarray*}
   & & \left.\sin\left(\frac{\pi{x}}{c}\frac{\partial}{\partial{z}}\right)(-\ln(1-e^z))\right|_{z=0}\,=\,
   -\left.\sin\left(Y\frac{\partial}{\partial{X}}\right)\ln{X}\right|_{z=0}\\
   &=& -\left.\textrm{arccot}\frac{X}{Y}\right|_{z=0}=-\left.\arctan\frac{Y}{X}\right|_{z=0}=-\arctan\frac{-\sin(\pi{x}/c)}{1-\cos(\pi{x}/c)}\\
   &=& \arctan\frac{\sin(\pi x/(2c))\cos(\pi x/(2c))}{\sin^2(\pi x/(2c))}.
\end{eqnarray*}
When $\sin(\pi{x}/(2c))\neq0$ or $x\neq0,\;x\neq2c$, the above expression can be written in the form:
\[\left.\sin\left(\frac{\pi{x}}{c}\frac{\partial}{\partial{z}}\right)(-\ln(1-e^z))\right|_{z=0}=
\arctan\cot\frac{\pi{x}}{2c}=\frac{\pi}{2}-\frac{\pi{x}}{2c}\qquad(0<x<2c).\]
Similarly, for $\sin(\pi{x}/(2c))>0$ (or $0<x<2c$), we have
\begin{eqnarray*}
   & & \left.\cos\left(\frac{\pi{x}}{c}\frac{\partial}{\partial{z}}\right)(-\ln(1-e^z))\right|_{z=0}\,=\,
   -\left.\cos\left(Y\frac{\partial}{\partial{X}}\right)\ln{X}\right|_{z=0}\\
   &=& -\left.\frac{1}{2}\ln(X^2+Y^2)\right|_{z=0}= -\frac{1}{2}\ln\left(\left(1-\cos\frac{\pi{x}}{c}\right)^2+\sin^2\frac{\pi{x}}{c}\right)=-\ln\left(2\sin\frac{\pi{x}}{2c}\right).
\end{eqnarray*}

\textbf{Lemma 3.3.} Let $m\in\mathbb{N}$. We have the following integral formula for $\ln{x}$:
\begin{equation}\label{37}
\underbrace{\int^x_0dx\cdots}_m\int^x_0\ln{x}\,dx=\frac{x^m}{m!}\left(\ln{x}-H_m\right),
\end{equation}
where $H_m$ denotes the familiar harmonic numbers defined by
\begin{equation}\label{har}
  H_m:=\sum^m_{j=1}\frac{1}{j}\qquad(m\in\mathbb{N}).
\end{equation}

\textbf{Proof.} Let $x\in\Omega\subset\mathbb{R}^1,\;n\in\mathbb{N}$. We can use the mathematical induction to prove the following formula
\begin{equation}\label{1*}
v\frac{d^nu}{dx^n}=\sum^n_{j=0}(-1)^j{n\choose{j}}\frac{d^{n-j}}{dx^{n-j}}\left(u\frac{d^jv}{dx^j}\right)\quad\forall{v,u}\in{C^n}(\Omega).
\end{equation}

Letting $u=x^n/n!$ in (\ref{1*}), we have
\begin{equation}\label{1*+}
v=\sum^n_{j=0}(-1)^j{n\choose{j}}\frac{d^{n-j}}{dx^{n-j}}\left(\frac{x^n}{n!}\frac{d^jv}{dx^j}\right)\quad\forall{v}\in{C^n}(\Omega).
\end{equation}
Taking $(\int_0^xdx)^n$ of both sides gives
\begin{equation}\label{1*++}
\underbrace{\int^x_0dx\cdots}_n\int^x_0vdx=\sum^n_{j=0}(-1)^j{n\choose{j}}\underbrace{\int^x_0dx\cdots}_j\int^x_0\frac{x^n}{n!}\frac{d^jv}{dx^j}dx
\quad\forall{v}\in{C^n}(\Omega).
\end{equation}
Letting $v=\ln{x}$ in (\ref{1*++}), then it is proved.

\textbf{Lemma 3.4.} Let $x\in\mathbb{R}^1$ with $0<x<\pi$. Then
\begin{equation}\label{37j}
\ln\frac{\sin x}{x}=-\sum^\infty_{k=1}\frac{1}{k}\frac{\zeta(2k)}{\pi^{2k}}x^{2k}.
\end{equation}

\textbf{Theorem 3.1.}  Let $m\in\mathbb{N},\,c>0$. The sum function $S_m(t)$ has the following recurrence property:
\begin{equation}\label{25}
\left.\cos\left(\frac{\pi{x}}{c}\frac{\partial}{\partial{z}}\right)S_m(e^z)\right|_{z=0}=S_m(1)
-\frac{\pi}{c}\int^x_0\left.\sin\left(\frac{\pi{x}}{c}\frac{\partial}{\partial{z}}\right)S_{m-1}(e^z)\right|_{z=0}dx,
\end{equation}
\begin{equation}\label{26}
\left.\sin\left(\frac{\pi{x}}{c}\frac{\partial}{\partial{z}}\right)S_m(e^z)\right|_{z=0}=
\frac{\pi}{c}\int^x_0\left.\cos\left(\frac{\pi{x}}{c}\frac{\partial}{\partial{z}}\right)S_{m-1}(e^z)\right|_{z=0}dx.
\end{equation}

\textbf{Proof.}  Taking $S(x)=S_{m-1}(x)/x$ in Theorem 2.5, then it is proved by using Lemma 3.1.

\textbf{Theorem 3.2.} For $m,r\in\mathbb{N}$, the sum function $S_m(t)$ has the following recurrence property:
\begin{eqnarray}\label{28}
                & & \left.\sin\left(\frac{\pi{x}}{c}\frac{\partial}{\partial{z}}\right)S_{m-1}(e^z)\right|_{z=0}\nonumber\\
                &=& \sum^{r-1}_{k=1}(-1)^{k-1}\frac{1}{(2k-1)!}\left(\frac{\pi{x}}{c}\right)^{2k-1}S_{m-2k}(1)\nonumber\\
                & & +\,(-1)^{r-1}\left(\frac{\pi}{c}\right)^{2r-1}
\underbrace{\int^x_0dx\cdots}_{2r-1}\int^x_0\left.\cos\left(\frac{\pi{x}}{c}\frac{\partial}{\partial{z}}\right)S_{m-2r}(e^z)\right|_{z=0}dx.
\end{eqnarray}
Here (and elsewhere in this work) an empty sum is to be interpreted (as usual) to be zero.

\textbf{Proof.}  We can use the mathematical induction to prove it. According to (\ref{26}) of Theorem 3.1, it is obviously tenable when $r=1$ in (\ref{28}). Now we inductively hypothesize that it is tenable when $r=K$. Using Theorem 3.1, then it is tenable when $r=K+1$, thus Theorem 3.2 is proved.

Similarly,

\textbf{Theorem 3.3.} For $m,r\in\mathbb{N}$, the sum function $S_m(t)$ has the following recurrence property:
\begin{eqnarray}\label{30}
                & & \left.\cos\left(\frac{\pi{x}}{c}\frac{\partial}{\partial{z}}\right)S_{m-1}(e^z)\right|_{z=0}\nonumber\\
                &=& \sum^{r-1}_{k=0}(-1)^k\frac{1}{(2k)!}\left(\frac{\pi{x}}{c}\right)^{2k}S_{m-2k-1}(1)\nonumber\\
                & & +\,(-1)^r\left(\frac{\pi}{c}\right)^{2r}
\underbrace{\int^x_0dx\cdots}_{2r}\int^x_0\left.\cos\left(\frac{\pi{x}}{c}\frac{\partial}{\partial{z}}\right)S_{m-2r-1}(e^z)\right|_{z=0}dx.
\end{eqnarray}

\textbf{Theorem 3.4.}  Let $r\in\mathbb{N}$, $x\in\mathbb{R}^1$ with $0<x<2c$. Then we have the following Fourier series expressions:
\begin{eqnarray}\label{36}
 \sum^\infty_{n=1}\frac{1}{n^{2r+1}}\cos\frac{n\pi{x}}{c}&=&
 \sum^{r-1}_{k=0}(-1)^{k}\frac{1}{(2k)!}\left(\frac{\pi{x}}{c}\right)^{2k}\zeta(2r+1-2k)\nonumber\\
 & & +\,(-1)^{r-1}\left(\frac{\pi}{c}\right)^{2r}\underbrace{\int^x_0dx\cdots}_{2r}\int^x_0\ln\left(2\sin\frac{\pi{x}}{2c}\right)\,dx
\end{eqnarray}
and
\begin{eqnarray}\label{36j}
 \sum^\infty_{n=1}\frac{1}{n^{2r}}\sin\frac{n\pi{x}}{c}&=&
 \sum^{r-1}_{k=1}(-1)^{k-1}\frac{1}{(2k-1)!}\left(\frac{\pi{x}}{c}\right)^{2k-1}\zeta(2r+1-2k)\nonumber\\
 & & +\,(-1)^r\left(\frac{\pi}{c}\right)^{2r-1}\underbrace{\int^x_0dx\cdots}_{2r-1}\int^x_0\ln\left(2\sin\frac{\pi{x}}{2c}\right)\,dx.
\end{eqnarray}

\textbf{Proof.} By using the relationship (\ref{23}), if $S_0(t)=-\ln(1-t)$, then $S_{2r-2k}(1)=\zeta(2r+1-2k)$. Letting $m=2r+1$ in Theorem 3.3, and by using Theorem 2.4, we have
\begin{eqnarray*}
 & & \sum^{\infty}_{n=1}\frac{1}{n^{2r+1}}\cos\frac{n\pi{x}}{c}\\
 &=& \sum^{r-1}_{k=0}(-1)^k\frac{1}{(2k)!}\left(\frac{\pi{x}}{c}\right)^{2k}\zeta(2r+1-2k)\\
 & & +\,(-1)^r\left(\frac{\pi}{c}\right)^{2r}
\underbrace{\int^x_0dx\cdots}_{2r}\int^x_0\left.\cos\left(\frac{\pi{x}}{c}\frac{\partial}{\partial{z}}\right)(-\ln(1-e^z))\right|_{z=0}dx.
\end{eqnarray*}
Then (\ref{36}) is proved using Lemma 3.2. Letting $m=2r,\;S_0(t)=-\ln(1-t)$ in Theorem 3.2, similarly we obtain (\ref{36j}).

In Theorem 3.4, since  $r\in\mathbb{N}$, the Fourier series expressions are also tenable at the endpoints of the interval $0<x<2c$.

\textbf{Theorem 3.5.} Let $r\in\mathbb{N}$, $x\in\mathbb{R}^1$ with $0\leq{x}\leq2c$. Then
\begin{eqnarray}\label{36'}
 \sum^\infty_{n=1}\frac{1}{n^{2r+1}}\cos\frac{n\pi{x}}{c}
 &=& \sum^{r-1}_{k=0}(-1)^k\frac{1}{(2k)!}\left(\frac{\pi{x}}{c}\right)^{2k}\zeta(2r+1-2k)\nonumber\\
 & & +\,(-1)^r\frac{1}{(2r)!}\left(\frac{\pi{x}}{c}\right)^{2r}\left(H_{2r}-\ln\frac{\pi{x}}{c}\right)\nonumber\\
 & & +\,(-1)^r2\sum^\infty_{k=1}\frac{(2k-1)!}{(2r+2k)!}\frac{\zeta(2k)}{(2\pi)^{2k}}\left(\frac{\pi{x}}{c}\right)^{2r+2k}
\end{eqnarray}
and
\begin{eqnarray}\label{yl7}
 \sum^\infty_{n=1}\frac{1}{n^{2r}}\sin\frac{n\pi{x}}{c}&=&
 \sum^{r-1}_{k=1}(-1)^{k-1}\frac{1}{(2k-1)!}\left(\frac{\pi{x}}{c}\right)^{2k-1}\zeta(2r+1-2k)\nonumber\\
 & & +\,(-1)^{r-1}\frac{1}{(2r-1)!}\left(\frac{\pi{x}}{c}\right)^{2r-1}\left(H_{2r-1}-\ln\frac{\pi{x}}{c}\right)\nonumber\\
 & & +\,(-1)^{r-1}2\sum^\infty_{k=1}\frac{(2k-1)!}{(2r+2k-1)!}\frac{\zeta(2k)}{(2\pi)^{2k}}\left(\frac{\pi x}{c}\right)^{2r+2k-1}.
\end{eqnarray}
Here the formula (\ref{yl7}) is a slightly modified version of a result proven in a significantly different way by Katsurada \cite[p. 81, Theorem 2]{Kats}.

\textbf{Proof.}  According to the expression (\ref{36}) of Theorem 3.4, using Lemma 3.3 and Lemma 3.4, if $r\in\mathbb{N}$, then in the closed interval $[0,2c]$ we have
\begin{eqnarray*}
   & & \sum^\infty_{n=1}\frac{1}{n^{2r+1}}\cos\frac{n\pi{x}}{c} \\
   &=& \sum^{r-1}_{k=0}(-1)^{k}\frac{1}{(2k)!}\left(\frac{\pi{x}}{c}\right)^{2k}\zeta(2r+1-2k)
   +\frac{(-1)^{r-1}}{(2r)!}\left(\frac{\pi{x}}{c}\right)^{2r}\ln2 \\
   & & +\,(-1)^{r-1}\left(\frac{\pi}{c}\right)^{2r}\underbrace{\int^x_0dx\cdots}_{2r}\int^x_0\left(\ln\frac{\pi{x}}{2c}
   +\ln\frac{\sin\frac{\pi{x}}{2c}}{\frac{\pi{x}}{2c}}\right)\,dx \\
   &=&\sum^{r-1}_{k=0}\frac{(-1)^{k}}{(2k)!}\left(\frac{\pi{x}}{c}\right)^{2k}\zeta(2r+1-2k)
   +\frac{(-1)^{r}}{(2r)!}\left(\frac{\pi{x}}{c}\right)^{2r}\left(H_{2r}-\ln\frac{\pi{x}}{c}\right)\\
   & &+\,(-1)^{r-1}\left(\frac{\pi}{c}\right)^{2r}\underbrace{\int^x_0dx\cdots}_{2r}\int^x_0
   -\sum^\infty_{k=1}\frac{1}{k}\frac{\zeta(2k)}{\pi^{2k}}\left(\frac{\pi{x}}{2c}\right)^{2k}\,dx\\
   &=& \sum^{r-1}_{k=0}\frac{(-1)^{k}}{(2k)!}\left(\frac{\pi{x}}{c}\right)^{2k}\zeta(2r+1-2k)
   +\frac{(-1)^{r}}{(2r)!}\left(\frac{\pi{x}}{c}\right)^{2r}\left(H_{2r}-\ln\frac{\pi{x}}{c}\right)\\
   & &+\,(-1)^r\sum^\infty_{k=1}\frac{1}{k}\frac{(2k)!\,\zeta(2k)}{(2\pi)^{2k}(2r+2k)!}\left(\frac{\pi x}{c}\right)^{2r+2k}.
\end{eqnarray*}
Thus (\ref{36'}) is proved. According to expression (\ref{36j}), similarly we have (\ref{yl7}) by using Lemma 3.3 and Lemma 3.4.

\section{Main results}
\noindent\setcounter{equation}{0}
\renewcommand\theequation{4.\arabic{equation}}

\textbf{Lemma 4.1.} (See Srivastava and Tsumura \cite[p. 329, Lemma 4]{sri}) For $\Re(s)>1$
\begin{equation}\label{yc2}
\sum^\infty_{n=1}\frac{1}{n^s}\cos\frac{2n\pi}{3}=\frac{1}{2}(3^{1-s}-1)\zeta(s);
\end{equation}
\begin{equation}\label{yc4}
\sum^\infty_{n=1}\frac{1}{n^s}\cos\frac{n\pi}{2}=2^{-s}(2^{1-s}-1)\zeta(s);
\end{equation}
\begin{equation}\label{yc1}
\sum^\infty_{n=1}\frac{1}{n^s}\cos\frac{n\pi}{3}=\frac{1}{2}(6^{1-s}-3^{1-s}-2^{1-s}+1)\zeta(s).
\end{equation}

\textbf{Lemma 4.2.} Let $r\in\mathbb{N}$, $x\in\mathbb{R}^1$ with $-2c\leq{x}\leq2c$ ($c>0$ be a given real number). Then we have the following Fourier series relationships:
\begin{eqnarray}\label{yl8}
 & & r\sum^\infty_{n=1}\frac{1}{n^{2r+1}}\cos\frac{n\pi{x}}{c}+\frac{\pi x}{2c}\sum^\infty_{n=1}\frac{1}{n^{2r}}\sin\frac{n\pi x}{c}\nonumber\\
 &=& \sum^{r-1}_{k=0}(-1)^k\,\frac{r-k}{(2k)!}\left(\frac{\pi{x}}{c}\right)^{2k}\zeta(2r+1-2k)\nonumber\\
 & & +\,(-1)^{r-1}\sum^\infty_{k=0}\frac{(2k)!}{(2r+2k)!}\frac{\zeta(2k)}{(2\pi)^{2k}}\left(\frac{\pi{x}}{c}\right)^{2r+2k},
\end{eqnarray}
which is a slightly modified version of a result proven in a significantly different way by Katsurada \cite[p. 81, Theorem 1]{Kats}.

\textbf{Proof.} By $\frac{\pi x}{c}\times(\ref{yl7})+2r\times(\ref{36'})$, since $\zeta(0)=-1/2$ and
\[\frac{(2k)!}{(2r+2k)!}=\frac{(2k-1)!}{(2r+2k-1)!}-2r\frac{(2k-1)!}{(2r+2k)!}\qquad(r,k\in\mathbb{N}),\]
then we can prove (\ref{yl8}).

By applying Lemma 4.2 with $x=c$, we have
\begin{eqnarray}\label{yl11}
 \zeta(2r+1)&=&\frac{2^{2r}}{r(2^{2r+1}-1)}\sum^{r-1}_{k=1}(-1)^{k-1}\frac{(r-k)\pi^{2k}}{(2k)!}\zeta(2r+1-2k)\nonumber\\
 & & +\,(-1)^r\frac{(2\pi)^{2r}}{r(2^{2r+1}-1)}\sum^\infty_{k=0}\frac{(2k)!}{(2r+2k)!}\frac{\zeta(2k)}{2^{2k}}\quad(r\in\mathbb{N}),
\end{eqnarray}
which was given earlier by Cvijovi\'{c} and Klinowski \cite[p. 1265, Theorem A]{Cvij}. In particular, when $r=1$, (\ref{yl11}) immediately yields Ewell's formula:
\[\zeta(3)=-\frac{4\pi^2}{7}\sum^\infty_{k=0}\frac{\zeta(2k)}{(2k+1)(2k+2)2^{2k}},\]
which was found by Euler in 1772 (See, e.g., \cite[p. 1080, Section 7]{Ay}), was rediscovered by Ramaswami \cite{Ra} and by Ewell \cite{Ew}.

Since
\[\zeta(2k)=(-1)^{k-1}\frac{(2\pi)^{2k}}{2(2k)!}B_{2k}\quad\mbox{and}\quad\zeta(-2k+1)=-\frac{1}{2k}B_{2k}\quad(k\in\mathbb{N}),\]
the last term on the right-hand side of (\ref{yl8}) can be written in the form:
\begin{eqnarray*}
  & & \sum^\infty_{k=r}(-1)^k\,\frac{r-k}{(2k)!}\left(\frac{\pi{x}}{c}\right)^{2k}\zeta(2r+1-2k) \\
  &=& \lim_{k\rightarrow0}a_k+\sum^\infty_{k=1}(-1)^{r+k}\,\frac{-k}{(2r+2k)!}\left(\frac{\pi{x}}{c}\right)^{2r+2k}\zeta(-2k+1) \\
  &=& \sum^\infty_{k=0}(-1)^{r+k}\,\frac{B_{2k}}{2(2r+2k)!}\left(\frac{\pi{x}}{c}\right)^{2r+2k} \\
  &=& (-1)^{r-1}\sum^\infty_{k=0}\frac{(2k)!}{(2r+2k)!}\frac{\zeta(2k)}{(2\pi)^{2k}}\left(\frac{\pi{x}}{c}\right)^{2r+2k},
\end{eqnarray*}
where
\[\lim_{k\rightarrow0}a_k=\lim_{k\rightarrow0}\left[(-1)^{r+k}\,\frac{-k}{(2r+2k)!}\left(\frac{\pi{x}}{c}\right)^{2r+2k}\zeta(-2k+1)\right]
=(-1)^r\,\frac{B_0}{2(2r)!}\left(\frac{\pi{x}}{c}\right)^{2r}.\]
Thus the Katsurada's formula (\ref{yl8}) and  Cvijovi\'{c}-Klinowski's formula (\ref{yl11}) can be written in the form:
\begin{eqnarray}\label{yl8g}
 & & r\sum^\infty_{n=1}\frac{1}{n^{2r+1}}\cos\frac{n\pi{x}}{c}+\frac{\pi x}{2c}\sum^\infty_{n=1}\frac{1}{n^{2r}}\sin\frac{n\pi x}{c}\nonumber\\
 &=& \sum^\infty_{k=0}(-1)^k\,\frac{r-k}{(2k)!}\left(\frac{\pi{x}}{c}\right)^{2k}\zeta(2r+1-2k)\quad(r\in\mathbb{N},\,|x|\leq2c);
\end{eqnarray}
\begin{equation}\label{yl11g}
 \zeta(2r+1)=\frac{2^{2r}}{r(2^{2r+1}-1)}\sum^\infty_{k=1}(-1)^{k-1}\frac{(r-k)\pi^{2k}}{(2k)!}\zeta(2r+1-2k)\quad(r\in\mathbb{N}).
\end{equation}

Since each series in (\ref{yl8}) is uniformly convergent with respect to $x$ on the interval $[-2c,2c]$ for $r\in\mathbb{N}$, multiplying the both side of (\ref{yl8}) by $(\pi x/c)^{2r-1}$, and executing termwise differentiation in them with respect to $x$, namely
$\partial_x[(\pi x/c)^{2r-1}\times(\ref{yl8})]$, then the sine series term is counteracted. Thus we can obtain the following theorem:

\textbf{Theorem 4.1.} Let $r\in\mathbb{N}$, $x\in\mathbb{R}^1$ with $-2c\leq{x}<2c$ (or $-2c\leq x\leq2c$ for $r\in\mathbb{N}\setminus\{1\}$). Then we have the following Fourier series relationship:
\begin{eqnarray}\label{yl15}
 & & r(2r-1)\sum^\infty_{n=1}\frac{1}{n^{2r+1}}\cos\frac{n\pi{x}}{c}+\frac{1}{2}\left(\frac{\pi x}{c}\right)^2\sum^\infty_{n=1}\frac{1}{n^{2r-1}}\cos\frac{n\pi x}{c}\nonumber\\
 &=& \sum^{r-1}_{k=0}(-1)^k\frac{(r-k)(2r+2k-1)}{(2k)!}\left(\frac{\pi{x}}{c}\right)^{2k}\zeta(2r+1-2k)\nonumber\\
 & & +\,(-1)^{r-1}\sum^\infty_{k=0}\frac{(2k)!(4r+2k-1)}{(2r+2k)!}\frac{\zeta(2k)}{(2\pi)^{2k}}\left(\frac{\pi{x}}{c}\right)^{2r+2k}.
\end{eqnarray}

In the special cases of (\ref{yl15}) when $r=1$, since
\begin{equation}\label{yl16}
  \sum^\infty_{n=1}\frac{1}{n}\cos\frac{n\pi{x}}{c}=-\ln\left(2\sin\frac{\pi x}{2c}\right)\qquad(0<x<2c),
\end{equation}
the Theorem 4.1 readily yield

\textbf{Theorem 4.2.} Let $x\in\mathbb{R}^1$ with $0<x<2c$. Then
\begin{eqnarray}\label{yl17}
 & & \sum^\infty_{n=1}\frac{1}{n^3}\cos\frac{n\pi{x}}{c}-\frac{1}{2}\left(\frac{\pi x}{c}\right)^2\ln\left(2\sin\frac{\pi x}{2c}\right)\nonumber\\
 &=& \zeta(3)+\sum^\infty_{k=0}\frac{2k+3}{(2k+1)(2k+2)}\frac{\zeta(2k)}{(2\pi)^{2k}}\left(\frac{\pi{x}}{c}\right)^{2k+2}.
\end{eqnarray}

\textbf{Theorem 4.3.} We have the following series representations for $\zeta(3)$
\begin{equation}\label{yl18}
  \zeta(3)=-\frac{\pi^2}{13}\ln3-\frac{4\pi^2}{13}\sum^\infty_{k=0}\frac{2k+3}{(2k+1)(2k+2)}\frac{\zeta(2k)}{3^{2k}};
\end{equation}
\begin{equation}\label{yl19}
  \zeta(3)=-\frac{2\pi^2}{35}\ln2-\frac{8\pi^2}{35}\sum^\infty_{k=0}\frac{2k+3}{(2k+1)(2k+2)}\frac{\zeta(2k)}{4^{2k}};
\end{equation}
\begin{equation}\label{yl20}
  \zeta(3)=-\frac{\pi^2}{6}\sum^\infty_{k=0}\frac{2k+3}{(2k+1)(2k+2)}\frac{\zeta(2k)}{6^{2k}}.
\end{equation}

\textbf{Proof.} Letting $x=2c/3$ in (\ref{yl17}), and by using (\ref{yc2}) we have (\ref{yl18}). Similarly, we can prove (\ref{yl19}) and (\ref{yl20}).

\textbf{Theorem 4.4.} For $r\in\mathbb{N}\setminus\{1\}$
\begin{eqnarray}\label{yl21}
 \zeta(2r+1)&=& -\frac{2\pi^2(3^{2r-2}-1)}{r(2r-1)(3^{2r+1}-1)}\,\zeta(2r-1)+\frac{2\times3^{2r}}{r(2r-1)(3^{2r+1}-1)}\nonumber\\
 & & \times\left[\,\sum^{r-1}_{k=1}(-1)^{k-1}\frac{(r-k)(2r+2k-1)}{(2k)!}\left(\frac{2\pi}{3}\right)^{2k}\zeta(2r+1-2k)\right.\nonumber\\
 & & \left.+\,(-1)^r\left(\frac{2\pi}{3}\right)^{2r}\sum^\infty_{k=0}\frac{(2k)!(4r+2k-1)}{(2r+2k)!}\frac{\zeta(2k)}{3^{2k}}\right].
\end{eqnarray}

\textbf{Proof.} Letting $x=2c/3$ in (\ref{yl15}), and by using (\ref{yc2}) we have (\ref{yl21}).

When $r=2$ and $3$ we obtain
\begin{equation}\label{yl24}
  \zeta(5)=\frac{41\pi^2}{363}\,\zeta(3)+\frac{8\pi^4}{363}\sum^\infty_{k=0}\frac{(2k+7)\zeta(2k)}{(2k+1)(2k+2)\cdots(2k+4)3^{2k}}.
\end{equation}
\begin{equation}\label{yl25}
\zeta(7)=\frac{2188\pi^2}{16395}\,\zeta(5)-\frac{18\pi^4}{5465}\,\zeta(3)
-\frac{64\pi^6}{16395}\sum^\infty_{k=0}\frac{(2k+11)\zeta(2k)}{(2k+1)(2k+2)\cdots(2k+6)3^{2k}}.
\end{equation}

\textbf{Theorem 4.5.} For $r\in\mathbb{N}\setminus\{1\}$
\begin{eqnarray}\label{yl22}
 \zeta(2r+1)&=& -\frac{(2^{2r-1}-2)\pi^2}{r(2r-1)(2^{4r+1}+2^{2r}-1)}\,\zeta(2r-1)+\frac{2^{4r+1}}{r(2r-1)(2^{4r+1}+2^{2r}-1)}\nonumber\\
 & & \times\left[\,\sum^{r-1}_{k=1}(-1)^{k-1}\frac{(r-k)(2r+2k-1)}{(2k)!}\left(\frac{\pi}{2}\right)^{2k}\zeta(2r+1-2k)\right.\nonumber\\
 & & \left.+\,(-1)^r\left(\frac{\pi}{2}\right)^{2r}\sum^\infty_{k=0}\frac{(2k)!(4r+2k-1)}{(2r+2k)!}\frac{\zeta(2k)}{4^{2k}}\right].
\end{eqnarray}

\textbf{Proof.} Letting $x=c/2$ in (\ref{yl15}), and by using (\ref{yc4}) we have (\ref{yl22}).

When $r=2$ and $3$ we obtain
\begin{equation}\label{yl26}
  \zeta(5)=\frac{157\pi^2}{1581}\,\zeta(3)+\frac{16\pi^4}{1581}\sum^\infty_{k=0}\frac{(2k+7)\zeta(2k)}{(2k+1)(2k+2)\cdots(2k+4)4^{2k}}.
\end{equation}
\begin{equation}\label{yl27}
\zeta(7)=\frac{14306\pi^2}{123825}\,\zeta(5)-\frac{64\pi^4}{41275}\,\zeta(3)
-\frac{128\pi^6}{123825}\sum^\infty_{k=0}\frac{(2k+11)\zeta(2k)}{(2k+1)(2k+2)\cdots(2k+6)4^{2k}}.
\end{equation}

\textbf{Theorem 4.6.} For $r\in\mathbb{N}\setminus\{1\}$
\begin{eqnarray}\label{yl23}
 \zeta(2r+1)&=& \frac{2\pi^2(6^{2r-2}-3^{2r-2}-2^{2r-2}+1)}{r(2r-1)(3^{2r}(2^{2r}+1)+2^{2r}-1)}\,\zeta(2r-1)\nonumber\\
 & & +\,\frac{2\times6^{2r}}{r(2r-1)(3^{2r}(2^{2r}+1)+2^{2r}-1)}\nonumber\\
 & & \times\left[\,\sum^{r-1}_{k=1}(-1)^{k-1}\frac{(r-k)(2r+2k-1)}{(2k)!}\left(\frac{\pi}{3}\right)^{2k}\zeta(2r+1-2k)\right.\nonumber\\
 & & \left.+\,(-1)^r\left(\frac{\pi}{3}\right)^{2r}\sum^\infty_{k=0}\frac{(2k)!(4r+2k-1)}{(2r+2k)!}\frac{\zeta(2k)}{6^{2k}}\right].
\end{eqnarray}

\textbf{Proof.} Letting $x=c/3$ in (\ref{yl15}), and by using (\ref{yc1}) we have (\ref{yl23}).

When $r=2$ and $3$ we obtain
\begin{equation}\label{yl28}
  \zeta(5)=\frac{8\pi^2}{87}\,\zeta(3)+\frac{\pi^4}{261}\sum^\infty_{k=0}\frac{(2k+7)\zeta(2k)}{(2k+1)(2k+2)\cdots(2k+4)6^{2k}}.
\end{equation}
\begin{equation}\label{yl29}
\zeta(7)=\frac{3124\pi^2}{29655}\,\zeta(5)-\frac{2\pi^4}{3295}\,\zeta(3)
-\frac{16\pi^6}{88965}\sum^\infty_{k=0}\frac{(2k+11)\zeta(2k)}{(2k+1)(2k+2)\cdots(2k+6)6^{2k}}.
\end{equation}

Since $\zeta(2k)\rightarrow1$ as $k\rightarrow\infty$, the general term in our series representations has the order estimate:
\[O(m^{-2k}\cdot k^{-2r+1})\qquad(k\rightarrow\infty;\quad m=3,4,6;\quad r\in\mathbb{N}),\]
whereas the general term in each of these earlier similar series representations has the order estimate:
\[O(2^{-2k}\cdot k^{-2r})\qquad(k\rightarrow\infty;\quad r\in\mathbb{N}).\]
Thus, even in the special case when $m=3$, the series representing $\zeta(2r+1)$ converges faster in (\ref{yl21}) than in (\ref{yl11}).

\textbf{Definition 4.1.} The Riemann Zeta function $\zeta(s)$ can be defined in the complete form (This definition is valid for all complex $s$):
\begin{equation}\label{rize}
 \zeta(s):=\frac{1}{s-1}\left(\frac{\partial}{\partial z}\right)^{1-s}\!\!\left.\frac{ze^z}{e^z-1}\right|_{z=0}\quad(s\in\Omega:=\{s\in\mathbb{C}|s\neq1\}).
\end{equation}
Here $\partial_z^{1-s},\,z\in\mathbb{R}^1$ is an abstract operators taking $\xi^{1-s}$ as the symbols, namely
\begin{equation}\label{cxsz}
 \left(\frac{\partial}{\partial z}\right)^{1-s}e^{\xi z}:=\xi^{1-s}e^{\xi z}\quad(\xi\in\mathbb{R}^1,\,s\in\mathbb{C}).
\end{equation}

In the special case when $\Re(s)>1$, by applying (\ref{e0}), (\ref{e1}) and (\ref{e1'}) to (\ref{rize}), we have
\begin{eqnarray*}
 \zeta(s)&=&\frac{1}{s-1}\left(\frac{\partial}{\partial z}\right)^{1-s}\!\!\left.\frac{ze^z}{e^z-1}\right|_{z=0}
 =\sum^\infty_{k=0}\left(\frac{\partial}{\partial z}\right)^{1-s}\!\!\left.\frac{ze^{(1+k)z}}{1-s}\right|_{z=0}\\
 &=&\sum^\infty_{k=0}\frac{\partial}{\partial\xi}e^{(1+k)\frac{\partial}{\partial\xi}}\left.\frac{\xi^{1-s}e^{\xi{z}}}{1-s}\right|_{\xi=0,z=0}\\
 &=&\sum^\infty_{k=0}\left(z+\frac{\partial}{\partial\xi}\right)e^{(1+k)\left(z+\frac{\partial}{\partial\xi}\right)}\left.\frac{\xi^{1-s}}{1-s}\right|_{\xi=0,z=0}\\
 &=&\sum^\infty_{k=0}\left(z+\frac{\partial}{\partial\xi}\right)\left.\frac{e^{(1+k)z}(\xi+1+k)^{1-s}}{1-s}\right|_{\xi=0,z=0}\\
 &=&\sum^\infty_{k=0}\left.\left(z\frac{e^{(1+k)z}(1+k)^{1-s}}{1-s}+e^{(1+k)z}(1+k)^{-s}\right)\right|_{z=0}\\
 &=&\sum^\infty_{k=0}\frac{1}{(1+k)^s}=\sum^\infty_{k=1}\frac{1}{k^s}\quad(\Re(s)>1).
\end{eqnarray*}
The result of this calculation shows that our definition is reasonable.

When $s=-n,\,n\in\mathbb{N}_0$ in (\ref{rize}), we obtain
\begin{equation}\label{hurn}
 \zeta(-n)=-\frac{1}{n+1}\frac{\partial^{n+1}}{\partial z^{n+1}}\!\left.\frac{ze^z}{e^z-1}\right|_{z=0}=-\frac{1}{n+1}B_{n+1}(1),
\end{equation}
where $B_n(x)$ are the Bernoulli polynomials defined by the generating functions:
\begin{equation}\label{bnld}
  \frac{ze^{xz}}{e^z-1}=\sum^\infty_{n=0}B_n(x)\frac{z^n}{n!}\quad(|z|<2\pi).
\end{equation}
It is easily seen from the definition (\ref{rize}) that
\begin{equation}\label{djd3}
 \lim_{s\rightarrow1}[(s-1)\zeta(s)]=\lim_{s\rightarrow1}\left(\frac{\partial}{\partial z}\right)^{1-s}\!\!\left.\frac{ze^z}{e^z-1}\right|_{z=0}
 =\left.\frac{ze^z}{e^z-1}\right|_{z=0}=B_0(1)=1.
\end{equation}
In other words, $\zeta(s)$ has a simple pole at $s=1$, and $B_0(1)=1$ is the residue of $\zeta(s)$ at the simple pole $s=1$.

\textbf{Definition 4.2.} Let $\zeta_z(s)$ be an analytic function defined by
\begin{equation}\label{rize'}
 \zeta_z(s):=\frac{1}{s-1}\left(\frac{\partial}{\partial z}\right)^{1-s}\!\!\frac{ze^z}{e^z-1}\quad(s\in\Omega:=\{s\in\mathbb{C}|s\neq1\},\,z\in\mathbb{R}^1).
\end{equation}

In the special case when $\Re(s)>2$ and $z\leq0$, we have
\begin{equation}\label{rize"}
 \zeta_z(s)=\frac{1}{1-s}\sum^\infty_{n=1}\frac{ze^{nz}}{n^{s-1}}+\sum^\infty_{n=1}\frac{e^{nz}}{n^s}\quad(\Re(s)>2,\,z\leq0).
\end{equation}

By Definition 4.1, we have $\left.\zeta_z(s)\right|_{z=0}=\zeta(s)$ and
\begin{equation}\label{wfze}
 \left.\frac{\partial^k}{\partial z^k}\zeta_z(s)\right|_{z=0}=\frac{s-1-k}{s-1}\zeta(s-k)\quad(k\in\mathbb{N},\,s\neq1).
\end{equation}

From the differential relation (\ref{wfze}), we have for $k\in\mathbb{N},\,s\neq1$
\[ \left.\frac{\partial^{2k-1}}{\partial x^{2k-1}}\cos\left(\frac{\pi x}{c}\frac{\partial}{\partial z}\right)\zeta_z(s)\right|_{z=0,x=0}=0;\]
\[ \left.\frac{\partial^{2k}}{\partial x^{2k}}\cos\left(\frac{\pi{x}}{c}\frac{\partial}{\partial{z}}\right)\zeta_z(s)\right|_{z=0,x=0}
=(-1)^k\left(\frac{\pi}{c}\right)^{2k}\frac{s-1-2k}{s-1}\zeta(s-2k).\]
Here $c>0$ is a given real number. Therefore, we have the following Taylor expansion in the neighborhood of $x=0$:
\begin{equation}\label{taze}
 \left.\cos\left(\frac{\pi{x}}{c}\frac{\partial}{\partial{z}}\right)\zeta_z(s)\right|_{z=0}=\sum^\infty_{k=0}\frac{(-1)^k}{(2k)!}\left(\frac{\pi{x}}{c}\right)^{2k}
 \frac{s-1-2k}{s-1}\zeta(s-2k)\quad(s\neq1).
\end{equation}
On the other hand, by using (\ref{rize"}) and Theorem 2.4, since
\begin{eqnarray*}
 & &\left.\cos\left(\frac{\pi{x}}{c}\frac{\partial}{\partial{z}}\right)(ze^{nz})\right|_{z=0}\\
 &=&\left.\cos\left(\frac{\pi{x}}{c}\frac{\partial}{\partial{z}}\right)z\cdot\cos\left(\frac{\pi{x}}{c}\frac{\partial}{\partial{z}}\right)e^{nz}\right|_{z=0}
-\left.\sin\left(\frac{\pi{x}}{c}\frac{\partial}{\partial{z}}\right)z\cdot\sin\left(\frac{\pi{x}}{c}\frac{\partial}{\partial{z}}\right)e^{nz}\right|_{z=0}\\
 &=&-\frac{\pi{x}}{c}\sin\frac{n\pi x}{c}\quad(n\in\mathbb{N}),
\end{eqnarray*}
we obtain the following summation formula for $\Re(s)>2,\,x\in\Omega\subset\mathbb{R}^1$:
\begin{equation}\label{sjqh}
 \left.\cos\left(\frac{\pi{x}}{c}\frac{\partial}{\partial{z}}\right)\zeta_z(s)\right|_{z=0}=\sum^\infty_{n=1}\frac{1}{n^s}\cos\frac{n\pi x}{c}
 +\frac{\pi{x}/c}{s-1}\sum^\infty_{n=1}\frac{1}{n^{s-1}}\sin\frac{n\pi x}{c}.
\end{equation}

By making use of (\ref{taze}), (\ref{sjqh}) and (\ref{yl8g}), we obtain the following theorem:

\textbf{Theorem 4.7.} Let $\Re(s)>1$. For $x\in\mathbb{R}^1$ with $|x|\leq2c$, we have
\begin{eqnarray}\label{zegh}
 & & s\sum^\infty_{n=1}\frac{1}{n^{s+1}}\cos\frac{n\pi x}{c}+\frac{\pi{x}}{c}\sum^\infty_{n=1}\frac{1}{n^s}\sin\frac{n\pi x}{c}\nonumber\\
 &=& \sum^\infty_{k=0}(-1)^k\frac{s-2k}{(2k)!}\left(\frac{\pi{x}}{c}\right)^{2k}\zeta(s+1-2k)\quad(\Re(s)>1),
\end{eqnarray}
which generalize the Katsurada's formula (\ref{yl8}) to the half complex plane $\Re(s)>1$; When $x=c$ in (\ref{zegh}), we obtain
\begin{equation}\label{tgck}
 \zeta(s+1)=\frac{2^s}{s(2^{s+1}-1)}\sum^\infty_{k=1}(-1)^{k-1}\frac{(s-2k)\pi^{2k}}{(2k)!}\zeta(s+1-2k)\quad(s\neq0),
\end{equation}
which generalize the Cvijovi\'{c}-Klinowski's formula (\ref{yl11}) to the whole complex plane.

Since each series in (\ref{zegh}) is uniformly convergent with respect to $x$ on the closed interval $[-2c,2c]$ for $\Re(s)>2$, multiplying the both side of (\ref{zegh}) by $(\pi x/c)^{s-1}$, and executing termwise differentiation in them with respect to $x$, namely $\partial_x[(\pi x/c)^{s-1}\times(\ref{zegh})]$, then the sine series term is counteracted. Thus we can generalize Theorem 4.1 to the half complex plane $\Re(s)>2$:

\textbf{Theorem 4.8.} Let $\Re(s)>2$. For $x\in\mathbb{R}^1$ with $|x|\leq2c$, we have
\begin{eqnarray}\label{zegh'}
 & & s(s-1)\sum^\infty_{n=1}\frac{1}{n^{s+1}}\cos\frac{n\pi{x}}{c}+\left(\frac{\pi{x}}{c}\right)^2\sum^\infty_{n=1}\frac{1}{n^{s-1}}\cos\frac{n\pi{x}}{c}\nonumber\\
 &=&\sum^\infty_{k=0}(-1)^k\frac{(s-2k)(s+2k-1)}{(2k)!}\left(\frac{\pi{x}}{c}\right)^{2k}\zeta(s+1-2k)\quad(\Re(s)>2).
\end{eqnarray}

When $x=2c/3,c/2,c/3$ in (\ref{zegh'}), by applying Lemma 4.1 to Theorem 4.8, we can generalize the Theorem 4.4, Theorem 4.5 and Theorem 4.6 to the whole complex plane:

\textbf{Theorem 4.9.} For $s\in\mathbb{C}$
\begin{eqnarray}\label{zegh1}
 & & s(s-1)\frac{3^{-s}-1}{2}\zeta(s+1)+\left(\frac{2\pi}{3}\right)^2\frac{3^{2-s}-1}{2}\zeta(s-1)\nonumber\\
 &=& \sum^\infty_{k=0}(-1)^k\frac{(s-2k)(s+2k-1)}{(2k)!}\left(\frac{2\pi}{3}\right)^{2k}\zeta(s+1-2k);
\end{eqnarray}
\begin{eqnarray}\label{zegh2}
 & & s(s-1)\frac{2^{-s}-1}{2^{s+1}}\zeta(s+1)+\left(\frac{\pi}{2}\right)^2\frac{2^{2-s}-1}{2^{s-1}}\zeta(s-1)\nonumber\\
 &=& \sum^\infty_{k=0}(-1)^k\frac{(s-2k)(s+2k-1)}{(2k)!}\left(\frac{\pi}{2}\right)^{2k}\zeta(s+1-2k);
\end{eqnarray}
\begin{eqnarray}\label{zegh3}
 & & s(s-1)\frac{6^{-s}-3^{-s}-2^{-s}+1}{2}\zeta(s+1)+\left(\frac{\pi}{3}\right)^2\frac{6^{2-s}-3^{2-s}-2^{2-s}+1}{2}\zeta(s-1)\nonumber\\
 &=& \sum^\infty_{k=0}(-1)^k\frac{(s-2k)(s+2k-1)}{(2k)!}\left(\frac{\pi}{3}\right)^{2k}\zeta(s+1-2k).
\end{eqnarray}

In other words, in the special case of Theorem 4.9 when $s=2r,\,r\in\mathbb{N}$, we can obtain the Theorem 4.4, Theorem 4.5 and Theorem 4.6.

Yan'an Second School, Yan'an 716000, Shaanxi, PR China

\emph{E-mail address}: guangqingbi@sohu.com

\end{CJK*}
\end{document}